\documentclass{siamltex}

\usepackage{graphicx}
\usepackage{amsmath,amsfonts}

\newcommand{\C}{\mathbb{C}}
\newcommand{\R}{\mathbb{R}}

\begin{document}

\title{Within-burst synchrony changes for coupled elliptic bursters}

\author{Abul Kalam al Azad\thanks{Mathematics Research Institute,
School of Engineering, Computing and Mathematics,
University of Exeter, Exeter EX4 4QF, U.K.}
 \and Peter Ashwin${}^*$ \\}

\maketitle

%%%%%%%%%%%%%%%%%%%%%%%%%%%%%%%%%%%%%%%%%%%%%%%%%%%%%%%%%%%%%%%%%%%%%%%%%%%%%%%%%%%%%%%%%%%%%%%%%%%%%%%%%%%%%%%%%%

\begin{abstract}
We study the appearance of a novel phenomenon for linearly coupled identical bursters: synchronized bursts where there are changes of spike synchrony within each burst.
The examples we study are for normal form elliptic bursters where there is a
periodic slow passage through a Bautin (codimension two degenerate Andronov-Hopf) bifurcation. This burster has a subcritical Andronov-Hopf bifurcation at the onset of repetitive spiking while end of burst occurs via a fold limit cycle bifurcation. We study synchronization behavior of two and three Bautin-type elliptic bursters for a linear direct coupling scheme. Burst synchronization is known to be prevalent behavior among such coupled bursters, while spike synchronization is more dependent on the details of the coupling.

We note that higher order terms in the normal form that do not affect the behavior of a single burster can be responsible for changes in synchrony pattern; more precisely, we find within-burst synchrony changes associated with a turning point in the spiking
frequency.
\end{abstract}

%%%%%%%%%%%%%%%%%%%%%%%%%%%%%%%%%%%%%%%%%%%%%%%%%%%%%%%%%%%%%%%%%%%%%%%%%%%%%%%%%%%%%%%%%%%%%%%%%%%%%%%%%%%%%%%%%%

~

\section{Introduction}

Elliptic bursting in a neuronal system is a recurrent alternation between active phases (large amplitude oscillations) and quiescent phases (small amplitude oscillations). This kind of rhythmic pattern can be found in rodent trigeminal neurons  \cite{RodentBurst}, thalamic relay and reticularis neurons \cite{ThalamicrelayBurst,ThalaReticularisBurst}, the primary afferent neurons in the brain stem circuits \cite{PrimaryAfferentBurst}, and neurons in many other areas of the brain. It is clearly of interest for neuronal population information encoding and transmission where several bursters fire within a population. Patterns of synchrony of elliptic bursters are may also be helpful in understanding firing patterns in more general types of burster \cite{ComBre,GolJosShi,Izhikevich,RinzelClassBurst}.

In a previous study of the synchronization of elliptic bursters, Izhikevich examined a pair of coupled ``normal form" elliptic bursters \cite{IzhiEllipticSyn} characterized by slow passage through a Bautin  (codimension two Andronov-Hopf) bifurcation. In that study, burst (slow activity pattern) synchronization between the bursters was found to be easily achievable, whereas spike (fast activity pattern) synchronization was harder to achieve. Other studies include
\cite{DroErm} who have examined nonlinearly coupled Bautin bifurcations though not in a bursting
setting and \cite{Izh00,ErmKop,SchDanSteBra} who have looked at various aspects of burst and spike synchronization for a variety of coupled burster models.

In this article we study spike synchronization for coupled Bautin-type elliptic bursters with more complicated phase (spiking) dynamics. It transpires that higher order terms that are not important in the normal form of a single burster can be responsible for nontrivial phase dynamics in coupled bursters. In particular, we observe and explain coexistence of and transitions between in-phase and anti-phase spiking within a single burst for two and more coupled bursters. This sheds light onto possible dynamical patterns of spike synchronization for coupled bursters in neuronal systems.

We discuss a normal form for coupled Bautin-type elliptic bursters and focus on burst and spike synchronization in a system of $n$ identical coupled bursters with $z_j\in\C$, $u_j\in\R$ and $j=1,\cdots,n$ given by
\begin{equation}\label{equation_main}
\left.\begin{aligned}
\dot{z}_{j} &= \left(u_{j}+i\omega\right)z_{j}+
B z_{j}|z_{j}|^{2}+ C z_{j}|z_{j}|^{4}+K_{j}\\
\dot{u}_{j} &= \eta (a-|z_{j}|)^{2},
\end{aligned}\right\}
\end{equation}
where $\omega,a,\eta\in\R$ and $B,C\in \C$ are fixed parameters, $K_{j}$ represents coupling. We set
\begin{equation}\label{equation_zetagamma}
B=2+i\zeta=2+i\frac{\sigma r_{m}^{2}}{2},~~C=-1+i\gamma=-1-i\frac{\sigma}{4}.
\end{equation}
We assume that the coupling term is
\begin{equation}
\label{equation_coupling}
K_j=(\kappa_{1}+i\kappa_{2})\sum_{k=1}^{n} c_{jk}z_{k},
\end{equation}
where $\kappa_1,\kappa_2\in\R$ are constant coupling parameters and $c_{jk}$ a constant connectivity matrix. For convenience here we take $c_{jk}=1$ for $j\neq k$, $c_{jj}=0$; i.e. all-to-all coupling. Biologically, although there are no rigorous reductions of specific bursters to this model, one can think of $z=x+iy$ as a fast variable $x$ that is analogous membrane voltage, $y$ that is analogous to the fast current, and a slow variable $u$ analogous to a slow adaptation current for a neuronal burster.

The article is structured as follows: in Section~\ref{sec_model} we discuss the
individual burster behavior for the model equations (\ref{equation_main}). In Section~\ref{sec_two_numerics} we discuss two such bursters ($n=2$), showing within-burst synchrony changes. These are analysed using a model system with assumed full burst synchrony in Section~\ref{sec_constrained} and slow-fast dynamics \cite{dissectionburst} to reduce to an equation for within-burst phase difference. Bifurcation analyses of this equation help one understand the observed dynamics of the full model. Section~\ref{sec_three_numerics} shows that we can observe similar dynamics in the three coupled elliptic burster system $n=3$ and conclude with a discussion of some dynamical and biological implications in Section~\ref{sec_discuss}.

\section{The model for coupled Bautin bursters}\label{sec_model}

Bursting is a multiple time scale phenomenon. In bursting, the fast dynamics of repetitive spiking is modulated by a slow dynamics of recurrent alternation between active and quiescent states. As explained in \cite{Izhikevich} one may obtain bursting from a variety of dynamical mechanisms; here we focus on elliptic bursters (\ref{equation_main}) with bursting behaviour by slow passage through a Bautin bifurcation; we briefly review the single burster dynamics.

\subsection{Normal form for Bautin bifurcation}

Suppose we have a Bautin bifurcation, namely a codimension two Andronov-Hopf bifurcation where
the criticality changes on varying an additional parameter. Then there is a normal form
that is locally topologically equivalent to the bifurcation, and this normal form
may be written \cite{KuznetsovBook} for $z=x+iy\in\C$ as
\begin{equation}
\dot{z}=Az+Bz|z|^{2}+Cz|z|^{4}+O(|z|^{6}),
\label{equation_1}
\end{equation}
where $A$, $B$, and $C$ are complex coefficients. If we write $A=A_r+iA_i$, $B=B_r+iB_i$, and $C=C_r+iC_i$ one can verify that an Andronov-Hopf bifurcation occurs as $A_r$ passes through $0$ and a change of criticality occurs where $B_r$ also passes through zero. The fourth order term is needed to determine the criticality at the degenerate point $B_r=0$. We write $B$ and $C$ as in (\ref{equation_zetagamma}). It can be shown that $\zeta,\gamma$ and $O(|z|^{6})$ terms do not affect the local branching dynamics of the system (\ref{equation_1}). We will however argue that $\zeta$ and $\gamma$ will influence the synchrony for two or more coupled elliptic bursters.

From (\ref{equation_1}) we obtain bursting dynamics \cite{Izh00,IzhiEllipticSyn,Izhikevich} by coupling the system to a slow variable $u\in\R$ that is the Andronov-Hopf parameter for the Bautin normal form, such that for $z$ small $u$ increases, while for $z$ large $u$ decreases:
\begin{equation}
\label{equation_2}
\left.\begin{aligned}
\dot{z} &= (u+i\omega)z+(2+i\zeta)z|z|^{2}+(-1+i\gamma)z|z|^{4}\\  
\dot{u} &= \eta (a-|z|^{2}).
\end{aligned}\right\}
\end{equation}
Note that $\eta \ll 1$ is the ratio of the fast to slow time scale. The system (\ref{equation_2}) exhibits bursting for $0<a<1$ while tonic spiking sets in for $a>1$. 

\begin{figure}[h]
\centering
\includegraphics[type=eps,ext=.eps,read=.eps,scale=0.45,angle=0,clip=]{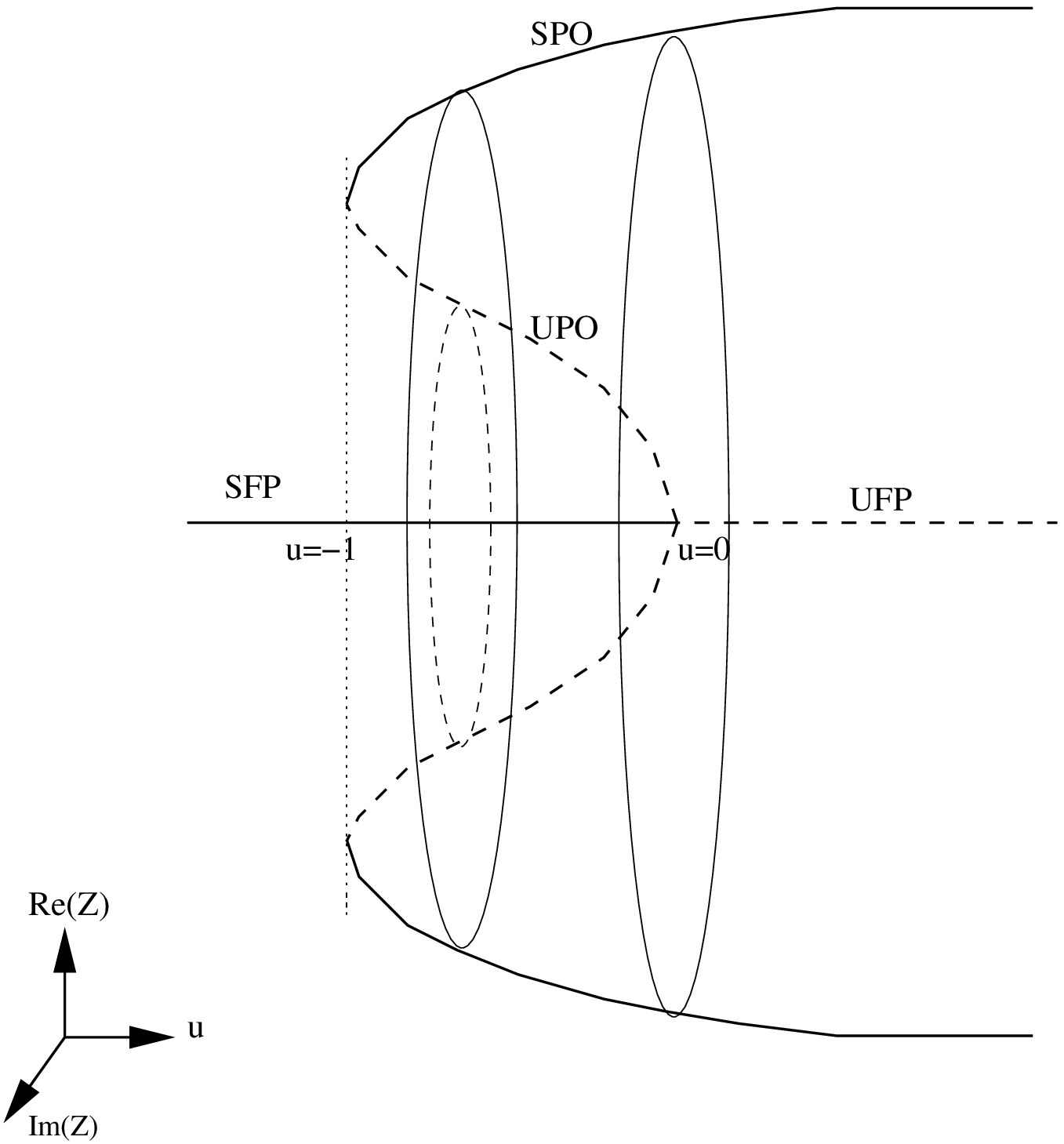}
\caption{\small{Schematic bifurcation diagram for $z$ in the system (\ref{equation_2}) on varying $u$. SFP denotes the stable fixed point, UFP unstable fixed point, SPO stable periodic orbit and UPO unstable periodic orbit. It is clearly seen that at $u=0$, the system undergoes subcritical Andronov-Hopf bifurcation, while saddle node bifurcation of limit cycles occur at $u=-1$.}}
\label{figure_1}
\end{figure}

In polar form, $z=r e^{i\theta}$, (\ref{equation_2}) becomes
\begin{equation}\label{equation_4}
\left.\begin{aligned}
\dot{r} &= u r+2r^{3}-r^{5}\\ 
\dot{\theta} &= \omega +\zeta r^2+\gamma r^4\\
\dot{u} &= \eta (a-r^{2}).
\end{aligned}\right\}
\end{equation}
In these coordinates it is clear that the fast subsystem undergoes an Andronov-Hopf bifurcation at $u=0$ and a limit cycle fold bifurcation (a saddle-node of limit cycles) at $u=-1$. At the saddle-node bifurcation of limit cycles, stable and unstable limit cycles coalesce. A bifurcation sketch for the system (\ref{equation_2}) is shown in the Figure~\ref{figure_1}. It is clear from this figure that the onset of periodic firing starts at a subcritical Andronov-Hopf bifurcation at $u=0$ with the emergence of a limit cycle. Likewise the steady state is reached via a saddle-node bifurcation of limit cycles at $u=-1$, where the stable limit cycle (solid line) meets the unstable limit cycle (dashed line) and eventually cancel each other at $u=-1$.

Note that during bursts, the limit cycles are {\em non-isochronous}; namely the interspike frequency $\dot{\theta}$ depends on $r$; there is a change in frequency of fast oscillation during the bursts. As this non-isochronicity does not affect the $r$ or $u$ dynamics, and hence the branching behaviour, it is not important for single bursters. The phase dynamics in (\ref{equation_4}) depends on amplitude $r$:
\begin{equation}
\dot{\theta}=\Omega(r)=\omega+\zeta r^{2}+\gamma r^{4}.
\label{equation_6}
\end{equation}

\begin{figure}[tp]
\centering
\includegraphics[type=eps,ext=.eps,read=.eps,scale=0.5]{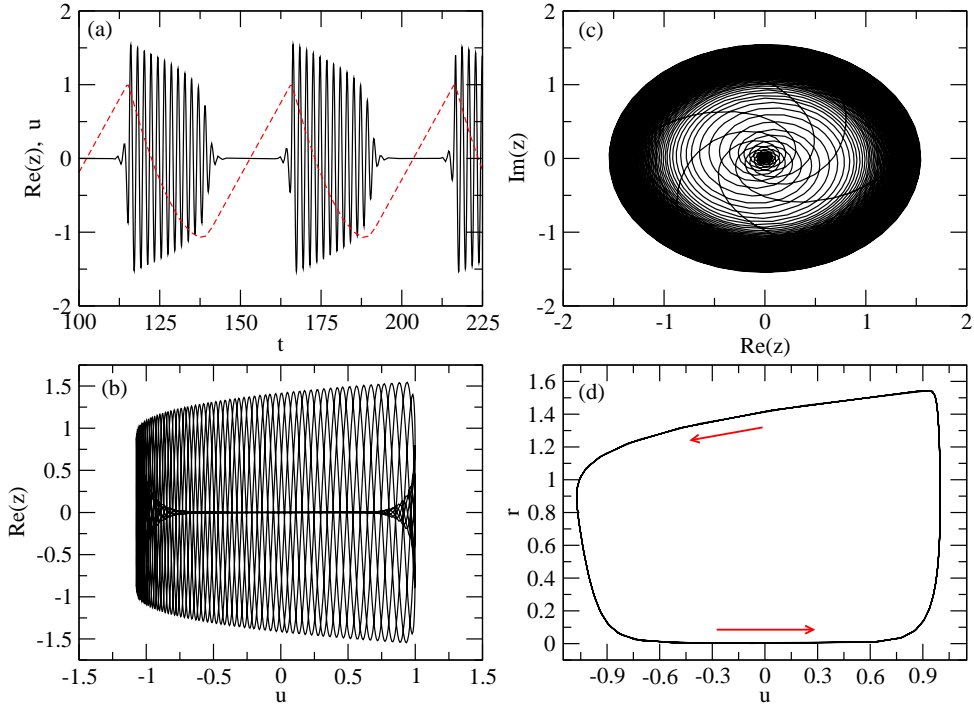}
\caption{\small{Dynamics of a single compartment Bautin Burster governed by ($\ref{equation_2}$) and ($\ref{equation_def_zetagamma}$). In panel (a), the timeseries of $Re(z)$ is shown with solid line and the corresponding slow variable, $u$, with dashed line.The parameters for the simulation are $\omega=3$, $a=0.8$, $\eta=0.1$, $\sigma=4$, and $r_{m}=1.35$. The arrows in (d) indicates the direction of change of the slow variable, $u$.}}
\label{figure_2}
\end{figure}

The non-trivial periodic orbits of the system (\ref{equation_4}) are $(r_{0},u_{0})=(\sqrt{a},a^{2}-2a)$ for $\eta=0$. Non-trivial periodic orbits, $r \neq 0$ correspond to periodic orbits of (\ref{equation_2}) with periodic spiking. The dynamics of the (\ref{equation_2}) is summarized in figure~\ref{figure_2} for parameters $\omega=3$, $\eta=0.1$, $a=0.8$, $\alpha=2$, $\beta=-1$, $\zeta=0$, and $\gamma=0$. Observe the slow passage effect \cite{slowpassage} apparent from Figure~\ref{figure_2}; although the stability calculation shows the Andronov-Hopf bifurcation occurs at $u=0$, but simulation shows a delayed bifurcation \cite{slowpassage}.

Note that we use parameters in (\ref{equation_zetagamma}) such that
\begin{equation}
\zeta=\frac{\sigma r_{m}^{2}}{2},~~\gamma=-\frac{\sigma}{4},
\label{equation_def_zetagamma}
\end{equation}

so that 
\begin{equation}
\frac{d\Omega}{dr} = \sigma r(r_{m}^{2}-r^{2}).
\label{equation_7}
\end{equation}
For the parameters ($\sigma, r_{m}$) it is clear that there is a turning point at $r=r_m$. The
parameter $\sigma$ can be interpreted as the magnitude of non-isochronicity for the phase dynamics while $r_{m}$ is a turning point for $\Omega$ on changing $r$.

\subsection{Coupled elliptic bursters}

We consider direct linear coupling of the normal form system (\ref{equation_2}) via the fast variables $z$ to give a coupled system of the form
(\ref{equation_main}, \ref{equation_coupling}) with coupling parameters $\kappa_{1}$ and $\kappa_{2}$. The coefficients $c_{jk}$ for the coupling term of (\ref{equation_main}) are the connectivity matrix; here we assume all-to-all coupling, namely 
$$
c_{jk}=\left\{
\begin{array}{ll}
1 &~~\mbox{ if }j\neq k\\
0 &~~\mbox{ otherwise }.
\end{array}\right.
$$

This form of coupling is analogous to the electrical (gap junction) coupling between synapses with phase shift expressed by the argument of $\kappa_{1}+i\kappa_{2}$. Positive $\kappa_{1}$ corresponds to excitatory coupling, while negative $\kappa_1$ corresponds to inhibitory coupling.

\section{Burst and spike synchronization for two coupled bursters}
\label{sec_two_numerics}

We numerically investigate the dynamics of a pair of coupled elliptic bursters governed by the system (\ref{equation_main}). Burst synchronization between the cells can be easily achieved for a wide range of parameter values with this system. In case of $\kappa_{2}=0$ and $\kappa_1>0$ (excitatory coupling), this generally generates inphase bursts, while antiphase bursts result from inhibitory coupling.\footnote{We write the system (\ref{equation_main}) using $z_{1}=x_{1}+i y_{1}$ and $z_{2}=x_{2}+i y_{2}$ for the purposes of numerical simulation. All the simulations were done with the interactive package XPPAUT~\cite{xppaut}. For integrations, the built-in adaptive Runge-Kutta integrator was used, and results were checked using the adaptive Dormand-Prince integrator.}

\begin{figure}[tbp]
\centering
\includegraphics[type=eps,ext=.eps,read=.eps,scale=0.15]{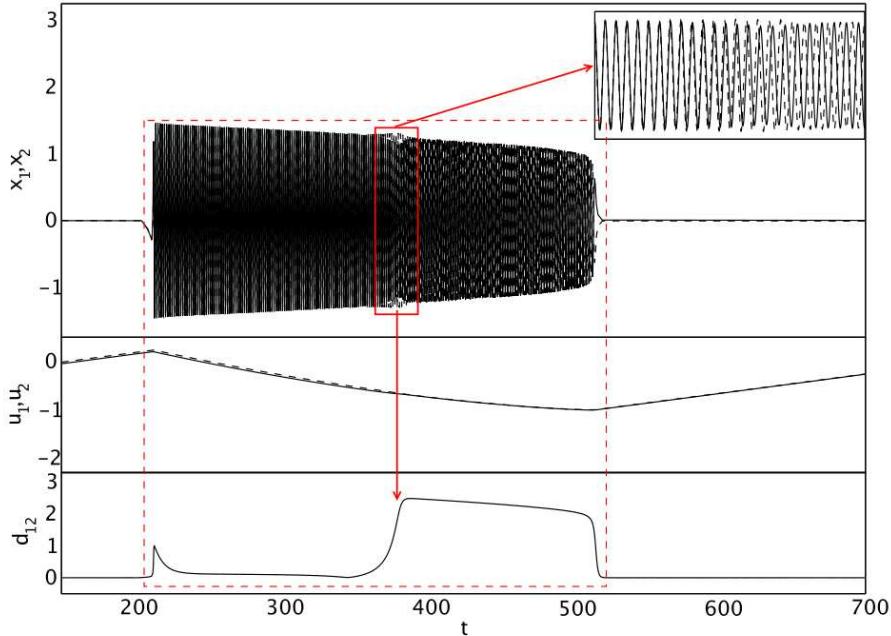}
\caption{\small{Within-burst synchrony change from stable inphase to stable antiphase states. This pattern repeats during each burst. This result is obtained from simulation of (\protect \ref{equation_main}), (\protect \ref{equation_zetagamma}) and (\protect \ref{equation_coupling}) for $n=2$, and parameters $\kappa_{1}=0.001$, $\kappa_{2}=0.2$, $\sigma=3$, $\eta=0.05$, $r_{m}=1.35$, $\omega=0.01$. Noise of amplitude  $10^{-5}$ was added to the fast subsystem. In this figure, the two coupled bursters are burst synchronized and the spikes become inphase at the beginning of the burst, but this synchrony pattern changes to antiphase near the middle of the burst. The inset in the topmost panel shows the region of the transition. Note that the initial transient and sudden change in the synchrony pattern along the burst profile are observable from $d_{12}$, where $d_{12}=0$ indicates inphase synchronization.}}
\label{figure_3}
\end{figure}

There is a ``within-burst synchrony change" observable within figure~\ref{figure_3}. The top panel shows $x_{1}$ and $x_{2}$. All transients were allowed to decay and the displayed pattern is repeated within each burst. A detail of the middle of the burst is shown in the top-right inset. The corresponding slow variables of the system, $u_{1}$ and $u_{2}$, are shown in the middle panel. The distinguishing solid and dashed traces correspond to the activity patterns of the two cells, respectively. The bottom shows the Euclidean distance 
$$
d_{12}=\sqrt{(x_{1}-x_{2})^{2}+(y_{1}-y_{2})^{2}+(u_{1}-u_{2})^{2}},
$$
between the two systems to show the presence ($d_{12}=0$) or absence ($d_{12}>0$) of synchrony. The values of the parameters used in the simulation are $\kappa_{1}=0.001$, $\kappa_{2}=0.2$, $\sigma=3$, $\eta=0.05$, $r_{m}=1.35$, $\omega=0.01$, so the two cells allowed to have same frequency. Wiener noise of amplitude $10^{-5}$ was added to the system so that the system does not stick in any unstable state.

The spikes are inphase at the beginning of the burst, but change to antiphase about the middle of the burst. The inset shows the region of this transition. This transition region may be shifted along the burst profile on changing $r_{m}$. Larger values of $r_{m}$ shift this transition leftward along the burst profile, and vice versa. This sudden change in the synchrony pattern along the burst profile is also captured by $d_{12}$.

\begin{figure}[tbp]
\centering
\includegraphics[type=eps,ext=.eps,read=.eps,scale=0.11]{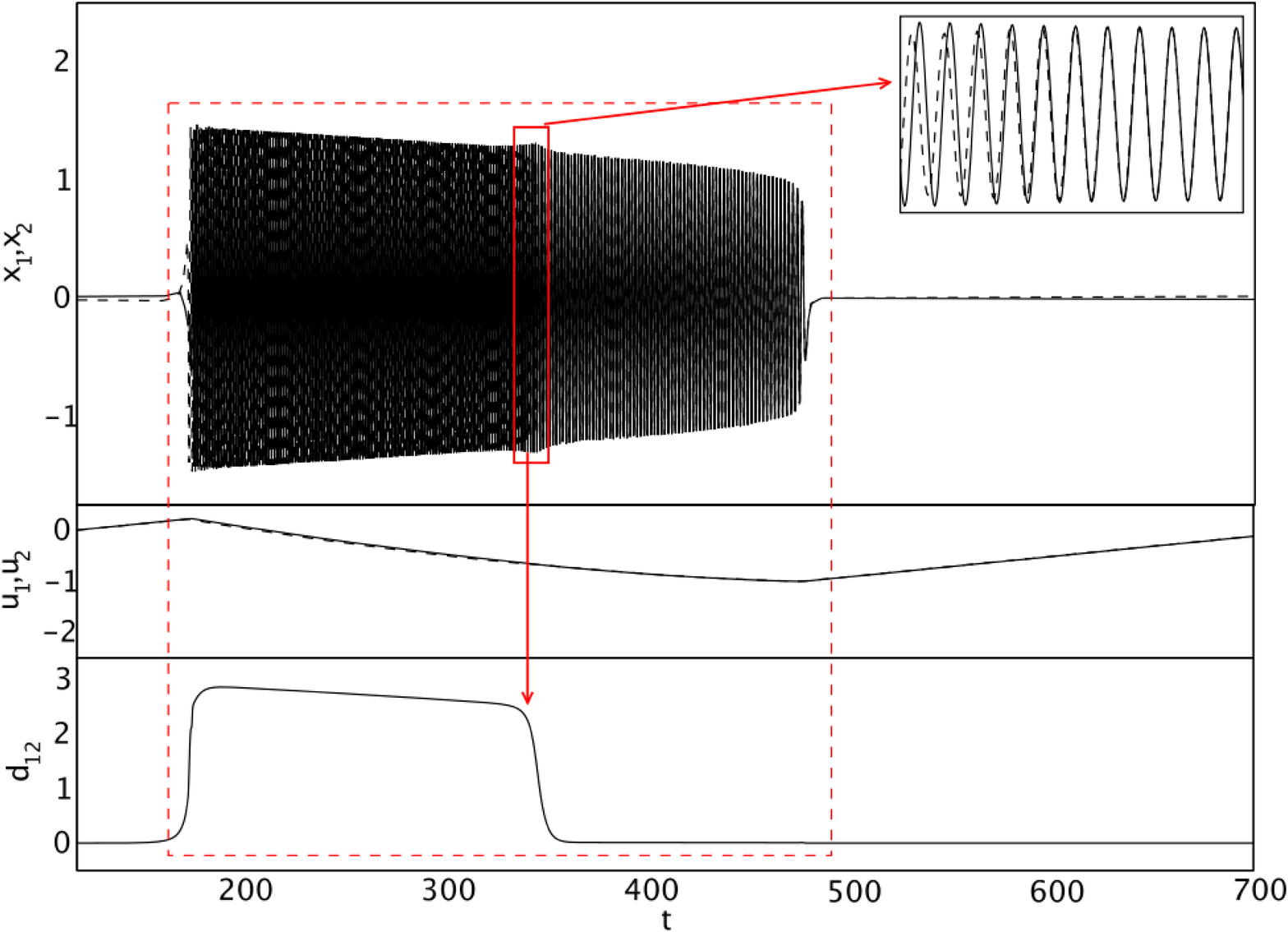}
\caption{\small{Within-burst synchrony change from stable antiphase to stable inphase states. The governing system and the parameters are as in figure~\protect \ref{figure_3} except $\kappa_{2}=-0.2$. The inset in the topmost panel shows the transition in detail. In the last panel, the bump in $d_{12}$ signifies the antiphase synchronization of the spikes within the synchronized burst. The corresponding slowly changing current variables, $u_{1}$ and $u_{2}$, are shown in the middle panel with solid and dotted lines, respectively.}}
\label{figure_4}
\end{figure}

We present another example of within-burst synchrony change for different parameter values in figure~\ref{figure_4}, where spikes of the two coupled cells start antiphase and change to inphase during the burst. The parameters for this are $\kappa_{1}=0.001$, $\kappa_{2}=-0.2$, $\eta=0.05$, $\sigma=3$, and $r_{m}=1.35$. As before, a low amplitude noise of order $10^{-5}$ was added to fast variables. The inset in the first panel shows the region of the transition. In the last panel, $d_{12}$ indicates that the burst is initially antiphase, and as it returns to $d_{12}=0$ there is a transition to inphase synchronization of the within-burst spikes of the two cells. The corresponding slowly changing current variables, $u_{1}$ and $u_{2}$, are shown in the middle panel. The overlapped solid and dashed lines imply the inphase burst synchronization of the coupled system.

\section{A burst-synchronized constrained model}
\label{sec_constrained}

Although it is possible to find within-burst synchrony changes within (\ref{equation_main}), we are not able to explain their existence analytically from the model. To overcome this,
we reduce the coupled system to a constrained problem where we assume burst synchronization, followed by a slow-fast decomposition. Using this we can explain how non-isochronicity and linear coupling can lead to within-burst synchrony changes.

\subsection{Two coupled bursters in polar coordinates}

Writing (\ref{equation_main}) with (\ref{equation_zetagamma}) and (\ref{equation_coupling}) in polar coordinates $z_j=r_j e^{i\theta_j}$ gives for $n=2$ the system
\begin{equation}
\label{equation_main_polar}
\left.\begin{aligned}
\dot{r}_{1} &= u_{1}r_{1}+2r^{3}_{1}-r^{5}_{1}+r_{2}(\kappa_{1} \cos(\theta_{2}-\theta_{1})-\kappa_{2} \sin(\theta_{2}-\theta_{1}))\\ 
\dot{\theta}_{1} &= \omega+\frac{1}{2}\sigma r_{m}^{2}r_{1}^{2}-\frac{1}{4}\sigma r_{1}^{4}+\frac{r_{2}}{r_{1}}(\kappa_{1}\sin(\theta_{2}-\theta_{1})+\kappa_{2}\cos(\theta_{2}-\theta_{1}))\\ 
\dot{u_{1}} &= \eta (a-r_{1}^{2})\\
\dot{r}_{2} &= u_{2}r_{2}+2r^{3}_{2}-r^{5}_{2}+r_{1}(\kappa_{1} \cos(\theta_{1}-\theta_{2})-\kappa_{2} \sin(\theta_{1}-\theta_{2}))\\ 
\dot{\theta}_{2} &= \omega+\frac{1}{2}\sigma r_{m}^{2}r_{2}^{2}-\frac{1}{4}\sigma r_{2}^{4}+\frac{r_{1}}{r_{2}}(\kappa_{1}\sin(\theta_{1}-\theta_{2})+\kappa_{2}\cos(\theta_{1}-\theta_{2}))\\ 
\dot{u_{2}} &= \eta (a-r_{2}^{2}).
\end{aligned}\right\}
\end{equation}  

We constrain the system to exact burst synchronization by setting:
\begin{equation}
\label{equation_constraint}
\left.\begin{aligned}
&u(t)=u_{1}(t)=u_{2}(t)\\
&\dot{u}=\dot{u}_1=\dot{u}_2 = \eta \left(a-\frac{1}{2}(|z_{1}|^{2}+|z_{2}|^{2}) \right).
\end{aligned}\right\}
\end{equation}

Thus, the system (\ref{equation_main_polar}) may be written with the constraint (\ref{equation_constraint}) and considering $\phi=\theta_{1}-\theta_{2}$ as
\begin{equation}
\label{equation_r1r2phi}
\left.\begin{aligned}
\dot{r}_{1} &= ur_{1}+2r^{3}_{1}-r^{5}_{1}+\kappa_{1}r_{2} \cos{\phi}+\kappa_{2}r_{2} \sin{\phi}\\ 
\dot{r}_{2} &= ur_{2}+2r^{3}_{2}-r^{5}_{2}+\kappa_{1}r_{1} \cos{\phi}-\kappa_{2}r_{1} \sin{\phi}\\ 
\dot{\phi} &= \frac{1}{2}\sigma r_{m}^{2}(r_{1}^{2}-r_{2}^{2})-\frac{1}{4}\sigma (r_{1}^{4}-r_{2}^{4})\\
&-\kappa_{1}\left(\frac{r_{1}^{2}+r_{2}^{2}}{r_{1}r_{2}}\right)\sin{\phi}-\kappa_{2}\left(\frac{r_{1}^{2}-r_{2}^{2}}{r_{1}r_{2}}\right)\cos{\phi}\\ 
\dot{u} &= \eta \left(a-\frac{1}{2}(r_{1}^{2}+r_{2}^{2})\right).
\end{aligned}\right\}
\end{equation}

As we are interested in synchrony changes, we define longitudinal and transverse coordinates
\begin{equation}
\left.\begin{aligned}
r_{l}&=(r_{1}+r_{2})/2\\
r_{t}&=(r_{1}-r_{2})/2.
\end{aligned}\right\}
\end{equation}
The system (\ref{equation_main_polar}) reduces to the four dimensional system
\begin{equation}
\label{equation_constrained}
\left.\begin{aligned}
\dot{r}_{l} &= ur_{l}+2r_{l}^{3}+6r_{l}r_{t}^{2}-r_{l}^{5}-10r_{l}^{3}r_{t}^{2}-5r_{l}r_{t}^{4}
\\
&+\kappa_{1}r_{l} \cos{\phi}-\kappa_{2}r_{t} \sin{\phi}\\ 
\dot{r}_{t} &= ur_{t}+6r_{l}^{2}r_{t}+2r_{t}^{3}-5r_{l}^{4}r_{t}-10r_{l}^{2}r_{t}^{3}-r_{t}^{5}
\\
&-\kappa_{1}r_{t} \cos{\phi}+\kappa_{2}r_{l} \sin{\phi}\\
\dot{\phi} &= 2\sigma r_{m}^{2}r_{l}r_{t}-2\sigma r_{l}r_{t}(r_{l}^{2}+r_{t}^{2})\\
&-2\kappa_{1}\frac{(r_{l}^{2}+r_{t}^{2})}{r_{l}^{2}-r_{t}^{2}}\sin{\phi}-4\kappa_{2}\frac{r_{l}r_{t}}{r_{l}^{2}-r_{t}^{2}}\cos{\phi}\\
\dot{u} &= \eta (a-(r_{l}^{2}+r_{t}^{2})).
\end{aligned}\right\}
\end{equation}

Here, $(r_{l},r_{t},\phi)$ govern the fast dynamics, and $u$ governs the slow dynamics. The system (\ref{equation_constrained}) is a reduced four-dimensional realization of the full system (\ref{equation_main_polar}) for a pair of coupled elliptic bursters. 

\subsection{Stability analysis of the burst constrained system}

In this section, we carry out a linear stability analysis of the fast sub-system of (\ref{equation_constrained}) about inphase and antiphase states with $r_{t}=0$ and $r_{l}=r$, which means both cells are burst synchronized with $r_{1}=r_{2}=r$. In the analysis, we assume the slow variable $u$ is a constant of the system by setting the time scale ratio, $\eta$, as a singularly perturbed parameter, i.e., $\eta=0$. The dynamics, as a result, is only governed by the fast spiking activity. For $r_{t}=0$ and $\eta=0$ we write $r_{l}=r$ as the stable nontrivial solution of (\ref{equation_constrained}) in the appropriate subspace
\begin{equation}
\label{equation_subspace}
\dot{r}_{l}=(u+2\kappa_{1} \cos{\phi})r_{l}+2r_{l}^{3}-r_{l}^{5}
\end{equation}
corresponding to bursting behaviour. Note that for small $|\kappa_{1}|$, this will have a solution close to the single burster case.

If we consider the fast system of (\ref{equation_constrained}) with $u$ between -1 and +1, then one can verify the existence of two solutions
\begin{itemize}
\item{Inphase} where $r_t=\phi=0$, $r_l=r$,
\item{Antiphase} where $r_t=0$, $\phi=\pi$, $r_l=r$,
\end{itemize}
where $r$ is the solution of 
\begin{equation}
\label{equation_solution_r}
u=r^{4}-2r^{2}-2\kappa_{1}\cos{\phi}.
\end{equation}

The Jacobian for the fast system at the inphase solution is block diagonal with a single real eigenvalue and a block
\begin{equation}
J_{in} = \left( \begin{array}{clcr}
u+6r^{2}-5r^{4}-\kappa_{1} & \kappa_{2}r \\
2\sigma r_{m}^{2}r-2\sigma r^{3}-\frac{4\kappa_{2}}{r} & -2\kappa_{1} 
\end{array} \right).
\label{equation_jin}
\end{equation}
Likewise, the Jacobian for the  fast system at the antiphase solution is also block diagonal with a single real eigenvalue and a block
\begin{equation}
J_{anti} = \left( \begin{array}{clcr}
u+6r^{2}-5r^{4}+\kappa_{1} & -\kappa_{2}r \\
2\sigma r_{m}^{2}r-2\sigma r^{3}+\frac{4\kappa_{2}}{r} & 2\kappa_{1} 
\end{array} \right).
\label{equation_janti}
\end{equation}
Note that the off-diagonal entries of the Jacobian matrices (\ref{equation_jin}) and (\ref{equation_janti}) depend on the imaginary part of the coupling coefficient, $\kappa_{2}$, and other system parameters. The real eigenvalues can be assumed negative because of stability of the solution of (\ref{equation_subspace}). 

The eigenvalues of the equation ($\ref{equation_jin}$) can be determined by examining the trace of the matrix ($\ref{equation_jin}$)
\begin{equation}
tr(J_{in})=u+6r^{2}-5r^{4}-3\kappa_{1}
\label{equation_trace_jin}
\end{equation}
and the determinant 
\begin{equation}
det(J_{in})=-2(u+6r^{2}-5r^{4})\kappa_{1}-2\sigma r^{2}(r_{m}^{2}-r^{2})\kappa_{2}+2\kappa_{1}^{2}+4\kappa_{2}^{2}.
\label{equation_det_jin}
\end{equation}
Many interesting insights into the inphase dynamics of ($\ref{equation_main}$), ($\ref{equation_zetagamma}$) and ($\ref{equation_coupling}$) may be extracted from ($\ref{equation_trace_jin}$) and ($\ref{equation_det_jin}$). Similarly, we can understand their antiphase counterparts from ($\ref{equation_janti}$) by examining
\begin{equation}
tr(J_{anti})=u+6r^{2}-5r^{4}+3\kappa_{1}
\label{equation_trace_janti}
\end{equation}
and
\begin{equation}
det(J_{anti})=2(u+6r^{2}-5r^{4})\kappa_{1}+2\sigma r^{2}(r_{m}^{2}-r^{2})\kappa_{2}+2\kappa_{1}^{2}+4\kappa_{2}^{2}.
\label{equation_det_janti}
\end{equation}

For simplicity, we consider the special case of weak coupling where $\kappa_{1}=0$, and $|\kappa_{2}| \ll 1$. In such a case, it may easily be seen that both $tr(J_{in})$ and $tr(J_{anti})$ in ($\ref{equation_trace_jin}$) and ($\ref{equation_trace_janti}$), respectively, are negative, as stability of the periodic solution of ($\ref{equation_subspace}$) means that $u+6r^{2}-5r^{4}<0$. So, from ($\ref{equation_trace_jin}$) and ($\ref{equation_trace_janti}$), $tr(J_{in})<0$, and $tr(J_{anti})<0$. For this weak coupling, it is also evident that $\big (tr(J_{in(anti)})\big )^{2}>4det(J_{in(anti)})$. Hence, the system will have a stable node for $det(J_{in(anti)})>0$ and a saddle for $det(J_{in(anti)})<0$.

To explain the within-burst synchrony change observed in figures $\ref{figure_3}$ and $\ref{figure_4}$, we write ($\ref{equation_det_jin}$) and ($\ref{equation_det_janti}$) to first order in $\kappa_{2}$, with $\kappa_{1}=0$ and small $\kappa_{2}>0$, as
\begin{equation}
det(J_{in})=-2\sigma r^{2}(r_{m}^{2}-r^{2})\kappa_{2}+O(\kappa_{2}^{2})
\label{equation_det_jin_reduced}
\end{equation}
and
\begin{equation}
det(J_{anti})=2\sigma r^{2}(r_{m}^{2}-r^{2})\kappa_{2}+O(\kappa_{2}^{2}).
\label{equation_det_janti_reduced}
\end{equation}

From equation ($\ref{equation_det_jin_reduced}$), if $r>r_{m}+O(\kappa_{2})$, then $det(J_{in})>0$. Together with the condition $tr(J_{in})<0$, this implies that the inphase solution is stable, whereas ($\ref{equation_det_janti_reduced}$) implies that the antiphase solution is unstable  for $r<r_{m}+O(\kappa_{2})$. We may derive approximate expressions for $r (=r_{l})$ where the bifurcations take place. We denote the bifurcation value for amplitude of the inphase solution by $r_{in}$, and the amplitude of the antiphase solution by $r_{anti}$. Note that $r_{in}$ may be obtained by equating $det(J_{in})$ to zero in the equation ($\ref{equation_det_jin}$) with $\kappa_{1}=0$ giving
\begin{equation}
det(J_{in})=-\kappa_{2}(2\sigma r_{in}^{2}-2\sigma r_{in}^{2}r_{m}^{2}-4\kappa_{2})=0.
\label{equation_ri_algebra}
\end{equation}
Now solving ($\ref{equation_ri_algebra}$),
\begin{equation}
r_{in}=r_{m}(1-\frac{\kappa_{2}}{\sigma r_{m}^{4}})+O(\kappa_{2}^{2}).
\label{equation_ri}
\end{equation}

Likewise, from equation ($\ref{equation_det_janti}$), the bifurcation point, $r_{anti}$, may be obtained as
\begin{equation}
r_{anti}=r_{m}(1+\frac{\kappa_{2}}{\sigma r_{m}^{4}})+O(\kappa_{2}^{2}).
\label{equation_ra}
\end{equation}
Note that $r$ depends on $u$ via ($\ref{equation_solution_r}$). So, the corresponding bifurcation points for in(anti)phase oscillations can be derived from ($\ref{equation_ri}$, $\ref{equation_ra}$) and ($\ref{equation_solution_r}$) as
\begin{equation}
u_{in}=r_{m}^{2}(r_{m}^{2}-2)+\frac{4\kappa_{2}}{\sigma r_{m}^{2}}(1-r_{m}^{2})+O(\kappa_{2}^{2}),
\label{equation_uin}
\end{equation} 
and
\begin{equation}
u_{anti}=r_{m}^{2}(r_{m}^2-2)-\frac{4\kappa_{2}}{\sigma r_{m}^{2}}(1-r_{m}^{2})+O(\kappa_{2}^{2}).
\label{equation_uanti}
\end{equation} 

\subsection{Synchrony bifurcations of the fast system}
\label{sec_two_bifs}

We now consider numerical bifurcation analyses of the fast system (\ref{equation_r1r2phi}) by taking $\eta$ as the singular perturbation parameter to take the fast system through single bursts and to compare with the asymptotic results found for $\kappa_{1}=0$.

\begin{figure}[tbp]
\centering
\includegraphics[type=eps,ext=.eps,read=.eps,scale=0.45]{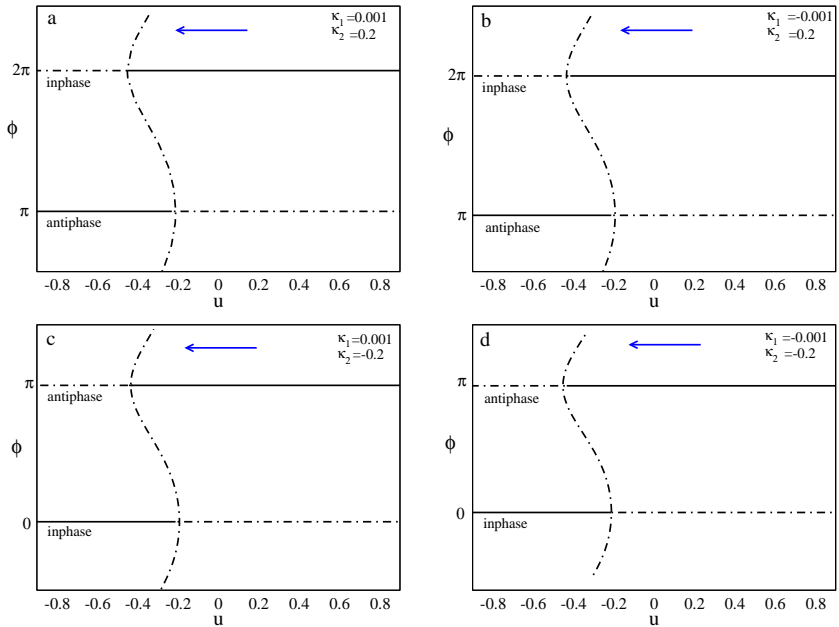}
\caption{\small{Bifurcation diagram of $\phi$ against $u$ for the burst-synchronized contrained system (\protect \ref{equation_r1r2phi}), where $u$ is a parameter that slowly decreases during each burst. The parameters are $\omega=3$, $\sigma=3$, $r_{m}=1.35$, and different $\kappa_{1}$ and $\kappa_{2}$ as indicated in the panels from (a) to (d). The solid lines represent stable solutions, while the unstable solutions are shown with dash-dotted lines. The arrow, running from right to left, shows the direction of the change of $u$ within a single burst.}}
\label{figure_5}
\end{figure}

To begin with we present in figure~\ref{figure_5} bifurcations of solutions of (\ref{equation_r1r2phi}) projected onto the phase difference, $\phi$, as $u$ is varied. The solid line represents the stable solutions, while the unstable solutions are shown with dash-dotted lines. The arrow, running from right to left, shows the direction of the change of $u$. What figure~\ref{figure_5}(a) shows is a burst that begins with stable inphase solution, and till almost half way through the burst, the inphase  solution remains stable and then the antiphase solutions gain stability. The coupling coefficients in this results are $\kappa_{1}=0.001$ and $\kappa_{2}=0.2$. The other parameter values are $\sigma=3$, $\omega=3$, and $r_{m}=1.35$. This behaviour agrees with the simulation result shown in the figure~\ref{figure_3} obtained from similar system parameters. Similarly, figure~\ref{figure_5}(c) explains what is found in the simulation in figure~\ref{figure_4}. Here, the spikes in the burst start off in stable antiphase and changes to stable inphase. Figure~\ref{figure_5}(b) and (d) show the results with $\kappa_{1}=-0.001$ but different $\kappa_{2}$. One interesting observation is the presence of the bistable region around the middle of the burst separating the stable inphase and antiphase solutions. This region occurs near the transition point ($r_{m}=1.35$) along the burst profile as predicted in the analysis in the previous section. These bifurcations show the robust coexistence of the inphase and antiphase synchrony patterns of the within burst spikes for a range of $u$, and within-burst synchrony changes of the coupled bursting system (\ref{equation_main}).

\begin{figure}[tbp]
\centering
\includegraphics[type=eps,ext=.eps,read=.eps,scale=0.35,clip=]{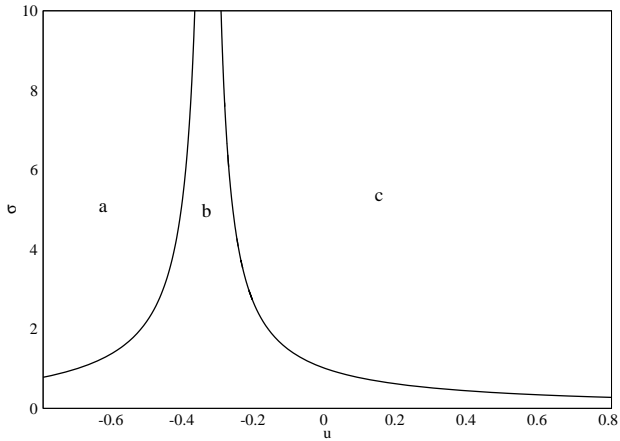}
\caption{\small{Two parameter bifurcation diagram of $\sigma$ against $u$ for the fast subsystem of (\protect \ref{equation_r1r2phi}). The other parameters are fixed at $\kappa_{1}=0.001$, $\kappa_{2}=0.2$, and $r_{m}=1.35$. There are stable inphase oscillations in region b,c, and stable antiphase oscillations in region a,b.}}
\label{figure_6}
\end{figure}

\begin{figure}[tbp]
\centering
\includegraphics[type=eps,ext=.eps,read=.eps,scale=0.35,clip=]{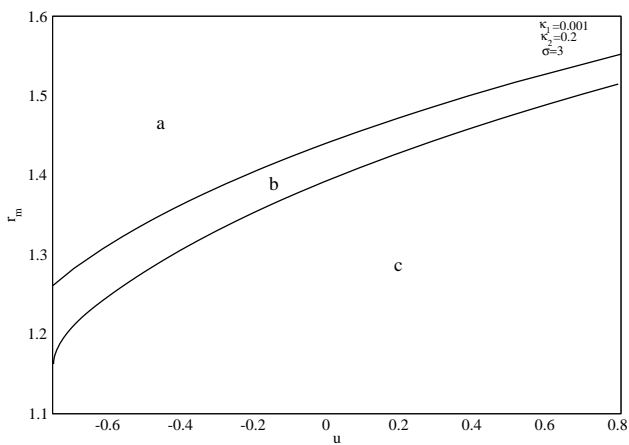}
\caption{\small{Two parameter bifurcation diagram of $r_{m}$ against $u$, for system and parameters as in figure $\protect \ref{figure_6}$ and $\sigma=3$. There are stable inphase oscillations in the region b,c, and stable antiphase oscillations in the region a,b.}}
\label{figure_7}
\end{figure}

Figure~\ref{figure_6} shows continuation in the $u\sigma$ plane. The range of  values of $u$ is within burst activity period. As $\sigma$ increases, the bistable region is seen to get narrower, in agreement with ($\ref{equation_ri}$, $\ref{equation_ra}$). Likewise, figure~\ref{figure_7} is obtained from parameters: $\kappa_{1}=0.001$, $\kappa_{2}=0.2$, and $\sigma=3$. What this result shows about the role of $r_{m}$ is that the position of the bistable region `b' may be shifted along the burst profile by varying $r_{m}$.

\begin{figure}[tbp]
\centering
\includegraphics[type=eps,ext=.eps,read=.eps,scale=0.35,clip=]{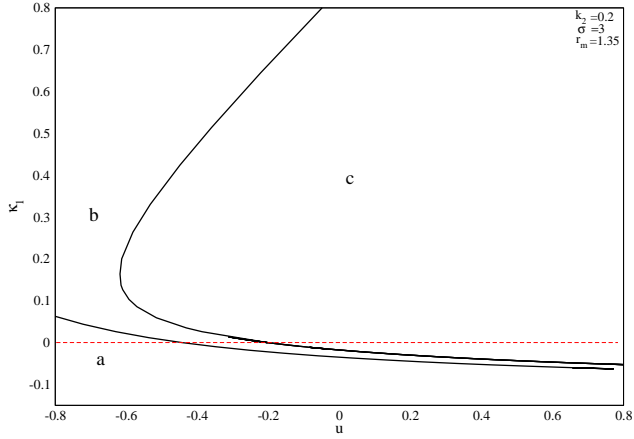}
\caption{\small{Bifurcation diagram of $\kappa_{1}$ against $u$ for system and parameters as in figure~$\ref{figure_5}$(a,b) with $\kappa_{2}=0.2$. There are stable inphase oscillations in the region b,c, and antiphase in the region a,b.}}
\label{figure_8}
\end{figure}

\begin{figure}[tbp]
\centering
\includegraphics[type=eps,ext=.eps,read=.eps,scale=0.35,clip=]{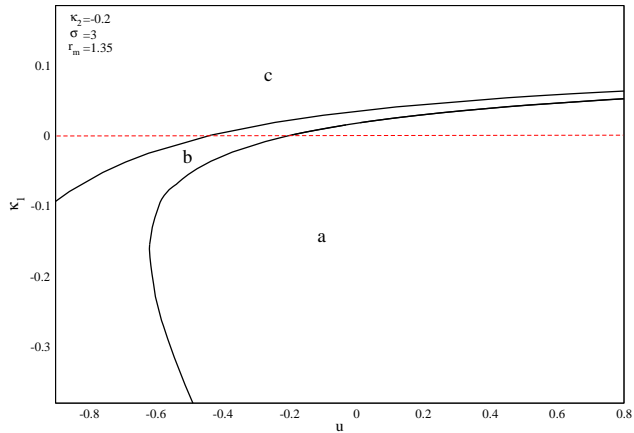}
\caption{\small{Two parameter bifurcation diagram in $\kappa_{1}$ against $u$ for system and parameters as in figure~$\ref{figure_5}$(c,d) with $\kappa_{1}$ varying and $\kappa_{2}=-0.2$. There are stable antiphase oscillations in the region a,b, and inphase in the region b,c.}}
\label{figure_9}
\end{figure}

The role of the coupling parameter $\kappa_{1}$ is shown in figures~\ref{figure_8} and \ref{figure_9} for two values of $\kappa_{2}$. The parameters in figure \ref{figure_8} are $\kappa_{2}=0.2$, $\sigma=3$, and $r_{m}=1.35$. For these parameter values, the behavior is similar as in figure~\ref{figure_5}(a) and (b). It is interesting to note that within-burst synchrony changes appear even for weak inhibitory values ($<0$) of $\kappa_{1}$. Moreover, stronger inhibitory values of $\kappa_{1}$ would mean only antiphase spike synchronization. Similarly, figure~\ref{figure_9} demonstrates similar dynamics to figure~\ref{figure_5}(c) and (d). Figure~\ref{figure_9} has parameters as those in figure~\ref{figure_8} except $\kappa_{2}=-0.2$. The excursion of the bistable region `b' above the dotted horizontal line indicates the appearance of within-burst synchrony changes for weak excitatory values, $\kappa_{1}>0$. Stronger  $\kappa_{1}$ results in inphase spike synchronization. 

Tables~\ref{table_comparison_1} and ~\ref{table_comparison_2} show the comparison of the inphase and antiphase bifurcation points, $r_{in}$, $u_{in}$ and $r_{anti}$, $u_{anti}$, for $\kappa_{1}=0$ and two values of $\kappa_{2}$ calculated from (\ref{equation_ri}, \ref{equation_uin}) and (\ref{equation_ra}, \ref{equation_uanti}), respectively, with those from simulations of systems (\ref{equation_r1r2phi}). Note that the bifurcation points obtained from the original system and those from the constrained system agree quite well.

\begin{table}[tbp]
\caption{Comparison of the bifurcation points, $r_{in}$, $r_{anti}$, $u_{in}$ and $u_{anti}$, obtained from simulations of system (\ref{equation_r1r2phi}) and those from equations (\ref{equation_ri}, \ref{equation_ra}, \ref{equation_uin}, \ref{equation_uanti}) for $\kappa_{1}=0$, $\kappa_{2}=0.2$, $\sigma=3$ and $r_{m}=1.35$.}
\centering
\begin{tabular}{|c|c|c|c|c|}
\hline
 & $r_{in}$ & $r_{anti}$ & $u_{in}$ & $u_{anti}$\\
\hline
\small{From system (\ref{equation_r1r2phi}) (figure~\ref{figure_8})} & 1.3210 & 1.376 & -0.4433 & -0.2027\\
\small{From equations (\ref{equation_ri}, \ref{equation_ra}, \ref{equation_uin}, \ref{equation_uanti})} & 1.3229 & 1.3771 & -0.4438 & -0.2032 \\
\hline
\end{tabular}
\label{table_comparison_1}
\end{table}

\begin{table}[tbp]
\caption{Comparison of the bifurcation points as in table~\ref{table_comparison_1} except $\kappa_{2}=-0.2$.} 
\centering
\begin{tabular}{|c|c|c|c|c|}
\hline
 & $r_{in}$ & $r_{anti}$ & $u_{in}$ & $u_{anti}$\\
\hline
\small{From system (\ref{equation_r1r2phi}) (figure~\ref{figure_9})} & 1.376 & 1.321 & -0.2027 & -0.4433\\
\small{From equations (\ref{equation_ri}, \ref{equation_ra}, \ref{equation_uin}, \ref{equation_uanti})} & 1.377 & 1.3229 & -0.2032 & -0.4438 \\
\hline
\end{tabular}
\label{table_comparison_2}
\end{table}

\begin{figure}[tbp]
\centering
\includegraphics[type=eps,ext=.eps,read=.eps,scale=0.35,clip=]{elliptic_k2vsu_k1positve_editted}
\caption{\small{Two parameter bifurcation diagram of $\kappa_{2}$ against $u$ for system and parameters as in figure~$\ref{figure_5}$(a,b) with $\kappa_{1}=0.001$ but $\kappa_{2}$ varied. There are stable inphase oscillations in the region b,c, and antiphase in the region a,b.}}
\label{figure_10}
\end{figure}

\begin{figure}[tbp]
\centering
\includegraphics[type=eps,ext=.eps,read=.eps,scale=0.35,clip=]{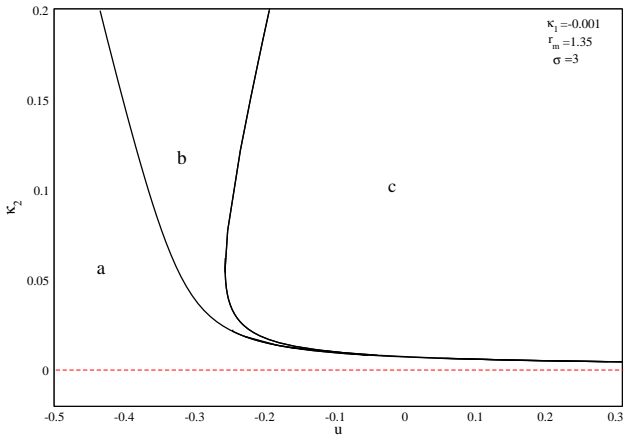}
\caption{\small{Bifurcation  diagram as in figure~$\ref{figure_10}$ but $\kappa_{1}=-0.001$ is varied. There are stable inphase oscillations in the region b,c, and antiphase in the region a,b.}}
\label{figure_11}
\end{figure}

Bifurcation diagrams in figures~\ref{figure_10} and \ref{figure_11} portray bifurcation diagram in $u\kappa_{2}$-space for two values $\kappa_{1}$. These figures show the role of $\kappa_{2}$ in spike synchronization. Figure~\ref{figure_10} uses the parameters: $\kappa_{1}=0.001$, $\sigma=3$, $r_{m}=1.35$. For negative and weak positive values of $\kappa_{2}$ (region c) only stable inphase solutions are present. Figure~\ref{figure_11} shows the bifurcation diagram for $\kappa_{1}=-0.001$. As before a, b is the region of stable antiphase, and b,c the stable inphase solutions. It may be observed that for negative and weak positive values of $\kappa_{2}$ (region a), one would see only stable antiphase synchronization of spikes in the burst. 

\begin{figure}[tbp]
\centering
\includegraphics[type=eps,ext=.eps,read=.eps,scale=0.18,clip=]{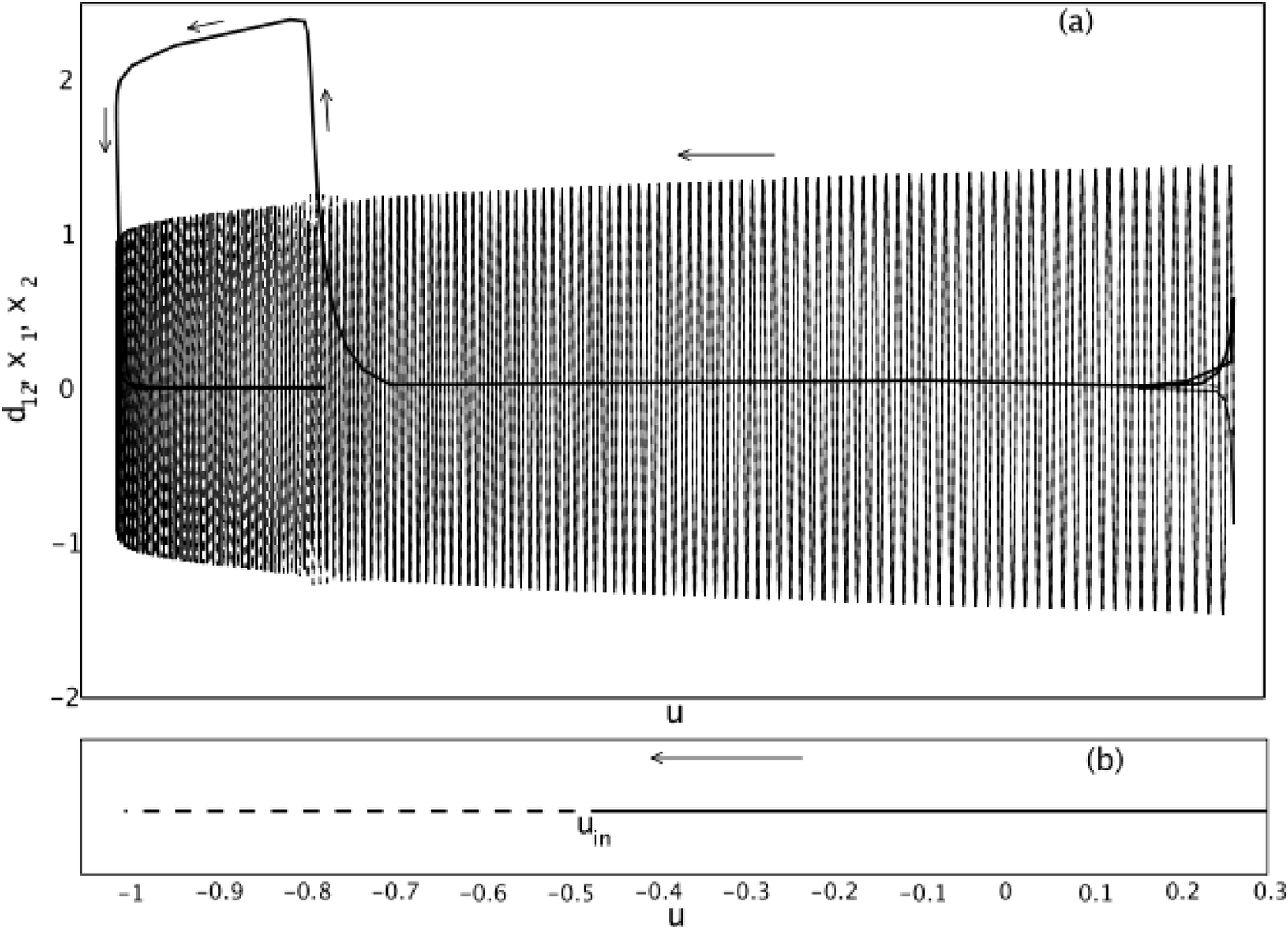}
\caption{\small{Comparison of the bifurcation point of within burst synchrony change (a) for the two coupled system (\ref{equation_main}, \ref{equation_zetagamma}, \ref{equation_coupling}) and (b) for the constrained system (\ref{equation_constraint}). In (a), the position of within burst synchrony change along $u$ can be observed from $d_{12}$ (thick solid line), where $u=(u_{1}+u_{2})/2$, is defined as such to agree with constrained slow variable $u$ in (\ref{equation_constrained}). The parameters for simulation in (a) are $\sigma=3$, $r_{m}=1.35$, $\kappa_{1}=0.001$, $\kappa_{2}=0.2$, $\eta=0.005$, $a=0.8$, and $\omega=0.0003$. A small amplitude noise of order $10^{-5}$ is added to the fast systems. Fast dynamics $x_{1}$ and $x_{2}$ are shown in solid and dashed lines, respectively. (b) depicts the inphase bifurcation solution branch of the figure~\ref{figure_5}(a). The solid line represents the stable inphase solution, and the unstable inphase solution is shown by dashed line. Note the difference in the inphase bifurcation points, $u_{in}$, and the loss of stability within burst.}}
\label{figure_uvsd12x1x2&phi_convertx}
\end{figure}

To end this section we present in figure~\ref{figure_uvsd12x1x2&phi_convertx} an interesting result comparing the bifurcation point of the within burst synchrony change between the original system (\ref{equation_main}, \ref{equation_zetagamma}, \ref{equation_coupling}) and the constrained system (\ref{equation_constraint}). Both systems have same coupling and system parameters. Note that the within burst synchrony change for the original system occurs at a more negative value of $u$ than that of the constrained system. This difference may be attributed to the slow passage effect of the fast system within the burst dynamics \cite{slowpassage}. We observed that the difference can be reduced by increasing noise amplitude added to the fast system. But eventually the large amplitude noise unsettles the stable solutions within the burst dynamics.

\section{Burst and spike synchronization for three bursters}
\label{sec_three_numerics}

\begin{figure}[tbp]
\centering
\includegraphics[type=eps,ext=.eps,read=.eps,scale=0.11]{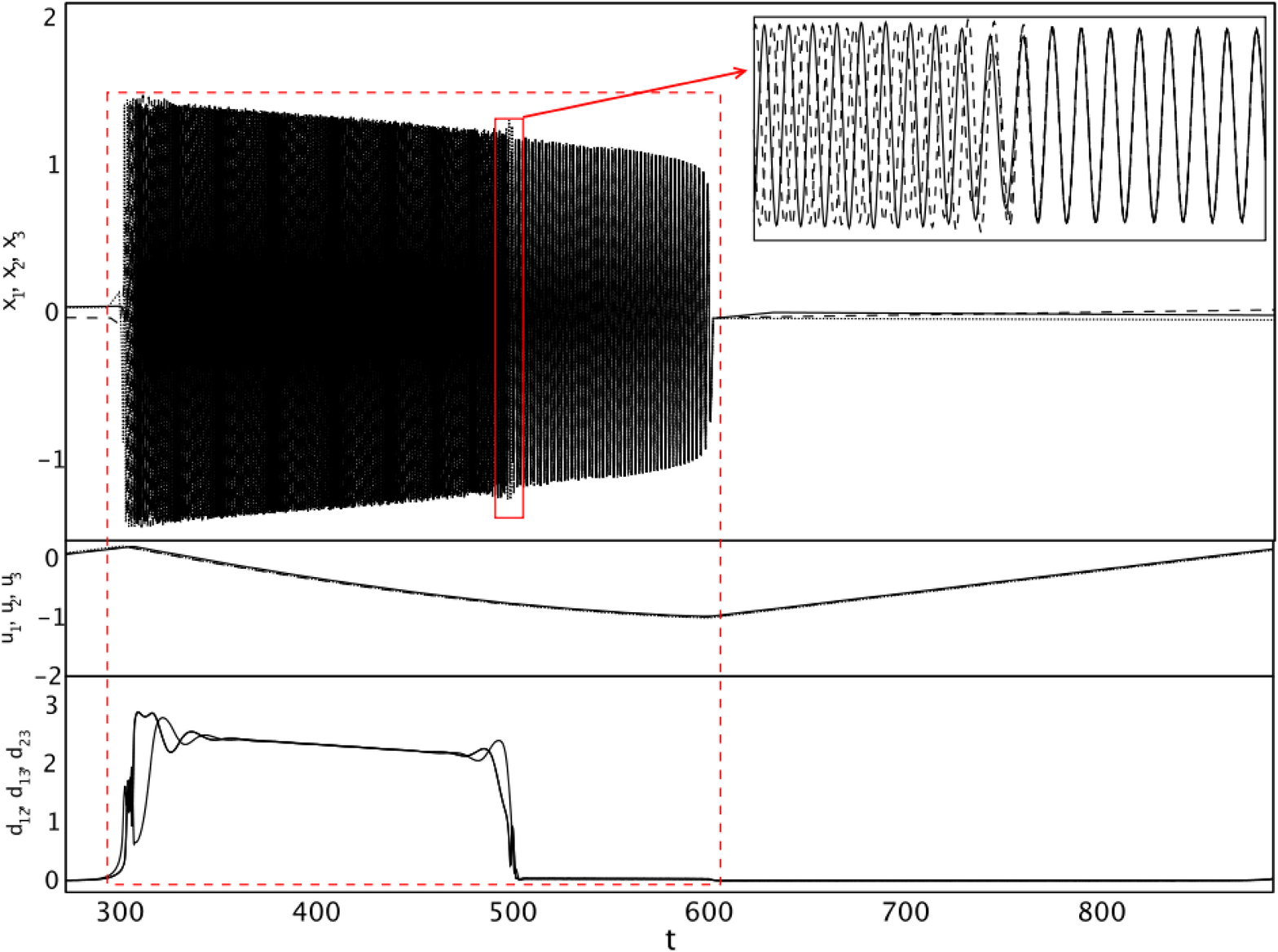}
\caption{\small{Within-burst synchrony change from stable antiphase to stable inphase states for three coupled Bautin-type elliptic bursters; see text for details. All three cells burst synchronously, but the fast spikes are antiphase at the beginning and inphase by the end of the burst.}}
\label{figure_12}
\end{figure}

We briefly demonstrate that within-burst synchrony changes are present in larger numbers of coupled bursters. In particular we look at three coupled Bautin-type elliptic bursters, i.e., (\ref{equation_main}), (\ref{equation_zetagamma}) and (\ref{equation_coupling}) with $n=3$.
A simulation of this system is shown in the figure~\ref{figure_12} with parameters $\omega=0.1$, $r_{m}=1.35$, $\sigma=5$, $\kappa_{1}=-0.001$, $\kappa_{2}=-0.2$, and additive noise of amplitude $10^{-5}$ to the fast variables.

The figure shows very similar behavior to the two-burster system with the difference that the
antiphase state is where the three bursters have a phase shift of $\frac{2\pi}{3}$ relative
to each other. The oscillations at the beginning of the burst are antiphase in this sense, and there is a transition to inphase during the burst. The inset in the top panel of the figure shows the transition in the spike synchrony pattern during the burst. The activity of the three different bursters is shown in solid, dashed and dash-dotted lines. The middle panel shows evolution of the corresponding slow variables, $u_{1}(t)$, $u_{2}(t)$ and $u_{3}(t)$. The third panel plots $d_{12}$, $d_{13}$ and $d_{23}$ that all must be zero for inphase synchronization, where

\begin{equation}
d_{ij}=\sqrt{(x_{i}-x_{j})^{2}+(y_{i}-y_{j})^{2}+(u_{i}-u_{j})^{2}},\nonumber\\
\end{equation}

with $i,j=1,2,3$. 

\section{Conclusion}
\label{sec_discuss}

We study the spiking dynamics of coupled elliptic bursters under direct linear coupling and find that within-burst synchrony changes are possible even for a simple normal form model, as long as terms that break isochronicity of the normal form are included. We observe that within-burst synchrony changes are stable and robust to changes in parameters. However, for identical bursters these within-burst changes are only easy to observe in the presence of noise; this is because in the absence of noise the system may become stuck in unstable synchronized states.

By reduction to fast-slow dynamics for the constrained burst-synchronised model we analyse the appearance of the within-burst synchrony change for two oscillators, and the influence of various system parameters. In particular we find that a turning point in the frequency $\Omega$
can be associated with the observed within-burst synchrony changes, analogous to bifurcations observed in systems of coupled weakly dissipative oscillators \cite{AshDan}. Moreover, we can find the approximate location of the transition between stable inphase and antiphase oscillations from bifurcation analysis of a reduced system. It will be a challenge to generalize this analysis to a case where the system is not burst-synchrony constrained. The examples we have illustrated in this paper are clearest for ``long" bursts where there are many oscillations during which the synchrony changes. Similar effects are presumably also present in shorter bursts, but are harder to observe because the changes in synchrony must occur over a small number of spikes to be observable. 

For larger populations of oscillators we expect there can be not just transitions between
inphase and antiphase during bursts, but also spontaneous changes in clustering, leading to robust but sensitive phase dynamics \cite{AshBor1,RabVarSelAba} and we believe this study gives some insight into the range of synchrony dynamics of coupled bursters in general. Better understanding of spike synchronization in more general coupled burster networks may lead to better understanding of potentially important new mechanisms for information processing and transmission by coupled neuronal bursters. This is discussed for example in~\cite{burstinformation2} where it is suggested that information transmission may occur via resonance between burst frequency and subthreshold oscillations.

\bibliographystyle{siam}
\bibliography{neuron_rev_new}

\end{document}